\documentclass[11pt,a4paper]{article}
\usepackage{amsmath}
\usepackage{amsthm}
\usepackage{amsfonts}
\parskip\smallskipamount
\theoremstyle{remark}
\newcommand{\E}{\mathbf{E}}

\newcommand{\Rp}{\mathbb{R}_+}
\newcommand{\Z}{\mathbb{Z}}
\newcommand{\Zp}{\mathbb{Z}_+}

\renewcommand{\AA}{\mathcal{A}}

\newcommand{\JJ}{\mathcal{J}}
\newcommand{\RR}{\mathcal{R}}
\newcommand{\NN}{\mathcal{N}}
\newcommand{\MM}{\mathcal{M}}

\newcommand{\gvn}{\mid}
\newcommand{\fc}[1]{\bar{#1}}
\renewcommand{\vec}[1]{\text{\boldmath $#1$}}

\parskip \smallskipamount

\title{Loss networks}
\author{Stan Zachary and Ilze Ziedins\\
  \textit{Heriot-Watt University and University of Auckland}
}
\begin{document}

\maketitle

\begin{quotation}\small\noindent
  We review the theory of loss networks, including recent results on
  their dynamical behaviour.  We give also some new results. 
  \vskip 0.1cm
  \noindent
%   \emph{Keywords:} loss networks.

%   \vskip 0.2cm

\end{quotation}

\section{Introduction}
\label{sec:introduction}

In a loss network calls, or customers, of various types are accepted
for service provided that this can commence immediately; otherwise
they are rejected.  An accepted call remains in the network for
some holding time, which is generally independent of the state
of the network, and throughout this time requires capacity
simultaneously from various network resources.  

The loss model was first introduced by Erlang as a model for the
behaviour of just a single telephone link (see Brockmeyer \emph{et al.} (1948).
The typical example remains that of a communications network, in which the
resources correspond to the links in the network, and a call of any type
requires, for the duration of its holding time, a fixed allocation of
capacity from each link over which it is routed (Kelly, 1986). This is
the case for
a traditional circuit-switched telephone network, but the model is
also appropriate to modern computer communications networks which
support streaming applications with minimum bandwidth
requirements (Kelly \emph{et al.}, 2000). 
There are also other examples: for instance, in a
cellular mobile network similar capacity constraints arise from the
need to avoid interference (Abdalla and Boucherie, 2002).

The mathematics of such networks has been widely studied, with
interest in both equilibrium and, more recently, dynamical behaviour.
Of particular importance are questions of call acceptance and capacity
allocation (for example, routing), with the aim of ensuring good
network performance which is additionally robust with respect to
variations in network parameters.  Call
arrival rates, in particular, may fluctuate greatly.  An excellent
review of the
state-of-the-art at the time of its publication is given by Kelly
(1991)---see also the many papers cited therein, and the later
survey by Ross (1995).

We take as our model the following.  Let $\RR$ denote the finite set
of possible call, or customer, types.  Calls of each type~$r\in\RR$
arrive at the network as a Poisson process with rate~$\nu_r$, and each
such call, if accepted by the network (see below), remains in it for a
\emph{holding time} which is exponentially distributed with
mean~$\mu_r^{-1}$.  We shall discuss later the extent to
which these assumptions, in particular the latter, are necessary.
Calls which are rejected do not retry and are simply considered lost.
All arrival processes and holding times are independent of one
another.  We denote the state of the network at time~$t$ by
$\vec{n}(t)=(n_r(t),\,r\in\RR)$, where $n_r(t)$ is the number of calls
of each type~$r$ in progress at that time.  The
process~$\vec{n}(\cdot)$ is thus Markov.  It takes values in some
state space $\NN\subset\Zp^R$, where $R=|\RR|$.  We assume $\NN$ to be
defined by a number of resource constraints
\begin{equation}
  \label{eq:1}
  \sum_{r\in\RR} A_{jr}n_r \le C_j, \qquad j\in \JJ,
\end{equation}
indexed in a finite set~$\JJ$, where the $A_{jr}$ and the $C_j$ are
nonnegative integers.  Typically we think of a call of each type~$r$
as having a simultaneous requirement, for the duration of its holding 
time, for $A_{jr}$ units of the capacity $C_j$ of
each resource~$j$; however, we show below that the resource
constraints~\eqref{eq:1} can also arise in other ways.  As noted
above, in
applications of this model to communications networks, the network
resources usually correspond to the \emph{links} in the network, and
when discussing the model in that context we shall generally find it
convenient to use this terminology.  We shall also
find it helpful to define, for each $r\in\RR$, the
parameter~$\kappa_r=\nu_r/\mu_r$; many quantities of interest depend
on $\nu_r$ and $\mu_r$ only through their ratio $\kappa_r$.

We shall say that a network is \emph{uncontrolled} whenever calls are
accepted subject only to the condition that the resulting state of the
network belongs to the set $\NN$.  Uncontrolled
networks are particularly amenable to mathematical analysis and are in
certain senses very well-behaved.  In addition, such a network has the
important \emph{insensitivity} property: the stationary
distribution of the process~$\vec{n}(\cdot)$ is unaffected by the
relaxation of the assumption that the call holding time distributions
are exponential, and depends on these holding time distributions only
through their
means.  This is essentially a consequence of the \emph{detailed
  balance} property considered in Section~\ref{sec:stat-dist}.

However, as we shall also see, the performance of uncontrolled
networks may be far from optimal.  A more general control strategy is
given by requiring that a call of type $r$, which arrives when the
state of the network (immediately prior to its arrival) is $\vec{n}$,
is accepted if and only if $\vec{n}\in\AA_r$ for some \emph{acceptance
  set}~$\AA_r$.  The sets~$\AA_r$ may be chosen so as to optimise, in
some appropriate sense, the network's performance.  Such networks do
not in general possess the insensitivity property described above.

Of interest in a loss network are both the stationary
distribution~$\pi$ and the dynamics of the process~$\vec{n}(\cdot)$.  For
the former it is usual to compute, for each $r$, the stationary
\emph{blocking probability}~$B_r$, that a call of type~$r$ is
rejected; here we shall find it slightly easier to work with the
stationary \emph{acceptance} (or \emph{passing})
\emph{probability}~$P_r=1-B_r$.  We note immediately that, by
Little's Theorem, the stationary expected number of
calls of each type~$r$ in the network is given by
\begin{equation}
  \label{eq:2}
  \E_\pi n_r = \kappa_r P_r,
\end{equation}
where $\kappa_r$ is as defined above.  Thus acceptance probabilities
may be regarded as one of the key performance measures in the
stationary regime.

It will be convenient to refer to the above model of a loss
network---in which arriving calls have fixed resource requirements and
in which the only control in the network is the ability to reject
calls---as the \emph{canonical model}.   When considering communications
networks, it is natural to extend this model by allowing also the
possibility of \emph{alternative routing}, in which calls 
choose their route according to the current state of the network.
Here the state space should properly be expanded to record the
number of calls of each type on each route (but see below).
% Where it is helpful, such a model may be incorporated into the
% framework already described by expanding the state space to record the
% number of calls of each type on each route, and by redefining arrival
% streams and acceptance sets appropriately.
We consider such models in Section~\ref{sec:mult-reso-models}.

In the case where we allow not only alternative routing, but also
\emph{repacking} of calls already in the network, the model simplifies
again, and it is once more only necessary for the state space to
record the number of calls of each type in progress.  Consider the
simple example of a communications network consisting of three links
with capacities~$C'_1$, $C'_2$, $C'_3$, and three call types, in which
calls of each type $r=1,2,3$ require \emph{either} one unit
of capacity from the corresponding link~$r$ \emph{or} one unit of
capacity from each of the other two links (in each case the
distribution of the call holding time is assumed to be the same).  If
repacking is allowed, the state of the system may be given by
$\vec{n}=(n_1,n_2,n_3)$ as usual, and it is easy to check that a call
of any type may be accepted if and only if the resulting state of the
network satisfies the constraints
\begin{displaymath}
  n_r + n_{r'} \le C'_r + C'_{r'}, \qquad r\neq r', \qquad r,r'\in\RR.
\end{displaymath}
This is therefore an instance of the uncontrolled network discussed
above, in which the coefficients~$A_{jr}$ and $C_j$ must be
appropriately defined.

Exact calculations for large loss networks typically exceed the
capabilities of even large computers, and we are thus led to consider
approximations.  Mathematical justification for these approximations
is usually
based on asymptotic results for one of two limiting schemes.  In the
first, which we shall refer to as the \emph{Kelly limiting scheme}
(see Kelly, 1986), the sets $\RR$, $\JJ$, the matrix $A = (A_{jr})$, and the
parameters~$\mu_r$ are held fixed, while the arrival rates~$\nu_r$ and
the capacities~$C_j$ are all allowed to increase in proportion to a
\emph{scale parameter}~$N$ which tends to infinity.  
In the second, which is known as the \emph{diverse routing limit} 
(see Whitt, 1985, and Ziedins and Kelly, 1989), the capacity of each 
resource is held constant, while the sets $\RR$ and $\JJ$ (and 
correspondingly the size of the matrix~$A$) are allowed to increase, and the 
arrival rates for call types requiring capacity at more than one 
resource to decrease, in such a way that the total traffic offered to each 
resource is also held constant (in particular, this requires that 
the arrival rate for any call type that requires capacity at more 
than one resource becomes negligible in the limit).
Results for the latter scheme in particular are used to justify
assumptions of independence in many approximations.

For each time~$t$, define $\vec{m}(t)=(m_j(t),\,j\in\JJ)$, where
$m_j(t)$ denotes the current occupancy, or usage, of each
resource~$j$ in a loss network.  Define also $\pi'$ to be the
stationary distribution of the process~$\vec{m}(\cdot)$.  In
particular, for the canonical model defined above, for each $t$,
\begin{equation}
  \label{eq:43}
  m_j(t) = \sum_{r\in\RR}A_{jr}n_r(t);
\end{equation}
here the process $\vec{m}(\cdot)$ takes values in the set
\begin{equation}
  \label{eq:46}
  \MM=\{\vec{m}\in\Zp^J\colon0\le{}m_j\le{}C_j,\;j\in\JJ\},
\end{equation}
where we
write $J=|\JJ|$,  and the distribution~$\pi'$ is given by
\begin{equation}
  \label{eq:44}
    \pi'(\vec{m}) = \sum_{\vec{n}\colon A\vec{n}=\vec{m}} \pi(\vec{n}),
  \qquad \vec{m}\in\MM.
\end{equation}
In general the process~$\vec{m}(\cdot)$ takes values in a space of
significantly lower dimension than that of the
process~$\vec{n}(\cdot)$.  This is
especially so in models of communications networks which incorporate
alternative routing.  It is a recurrent theme in the
study of loss networks that, in general, at least approximately
optimal control of a network is obtained by basing admission decisions
and, in communications networks, routing decisions, solely on the
state of the process~$\vec{m}(\cdot)$ at the arrival time of each
call.  Further, in this case, a knowledge of the distribution~$\pi'$
is sufficient to determine call acceptance probabilities.  We shall
also see that good estimates of $\pi'$ are generally given by assuming
its (approximate) factorisation as
\begin{equation}
  \label{eq:45}
    \pi'(\vec{m}) = \prod_{j\in\JJ} \pi'_j(m_j),
\end{equation}
where each $\pi'_j$ is normalised to be a probability distribution.
This is a further recurrent theme in the study of loss networks.

In Section~\ref{sec:rev} we consider the stationary behaviour of
uncontrolled networks, reviewing both exact results and approximations
for large networks.  Our approach is based on the use of an
elegant recursion due to Kaufman (1981) and to Dziong and Roberts
(1987) which delivers all the classical results in regard to, for
example, stationary acceptance probabilities, with a certain
simplicity.

More general networks are studied in
Section~\ref{sec:contr-loss-netw}.  In Section~\ref{sec:simple-model}
we study the problem of optimal control in a single-resource network,
where a reasonably tractable analysis of stationary behaviour is again
possible, and where we show that either exactly or approximately
optimal control may be obtained with the use of strategies based on
\emph{reservation} parameters.  In Section~\ref{sec:mult-reso-models}
we consider multiple-resource networks, allowing in particular the
possibility of alternative routing.  We again derive approximations
which are known to work extremely well in practice.  In
Section~\ref{sec:dyn} we consider the dynamical behaviour of large
loss networks.  This is important for the study of the long-run, and
hence also the equilibrium, behaviour of networks in the case where a
direct equilibrium analysis is impossible.  The study of network
dynamics is also the key to understanding their stability.   Finally, in
Section~\ref{sec:recent-devel-vari} we mention some wider models and
discuss some open problems.

\section{Uncontrolled loss networks: stationary behaviour}
\label{sec:rev}

We study here the stationary behaviour of the uncontrolled network
introduced above, in which calls of any type are accepted subject only
to the condition that the resulting state~$\vec{n}$ of the network
belongs to the state space~$\NN$ defined by the capacity
constraints~\eqref{eq:1}.  In particular we shall see, in
Section~\ref{sec:kdr-recursion} and subsequently, that most quantities
of interest, in particular acceptance probabilities, may be
calculated, exactly or approximately,
without the need to calculate the full stationary distribution~$\pi$
of the process~$\vec{n}(\cdot)$.

\subsection{The stationary distribution}
\label{sec:stat-dist}

For each $r\in\RR$, let $\vec{\delta}_r$ be the vector whose $r$th
component is $1$ and whose other components are $0$.  Recall that,
under the assumptions introduced above, $\vec{n}(\cdot)$ is a Markov
process.  For $\vec{n},\vec{n}-\vec{\delta}_r\in\NN$ and $r\in\RR$,
its transition rates
between $\vec{n}$ and $\vec{n}-\vec{\delta}_r$ are $n_r\mu_r$ and
$\nu_r$.  It thus follows that the stationary distribution~$\pi$ of
the process~$\vec{n}(\cdot)$ is given by the solution of the
\emph{detailed balance equations}
\begin{equation}
  \label{eq:3}
  \pi(\vec{n})n_r\mu_r = \pi(\vec{n}-\vec{\delta}_r)\nu_r,
  \qquad r \in \RR, \quad \vec{n}\in\NN,
\end{equation}
where, here and elsewhere, we make the obvious convention that
$\pi(\vec{n}-\vec{\delta}_r)=0$ whenever $n_r=0$.  That is, 
\begin{equation}
  \label{eq:4}
  \pi(\vec{n})
  = G^{-1}\prod_{r\in\RR}\frac{\kappa_r^{n_r}}{n_r!},
  \qquad \vec{n}\in\NN,
\end{equation}
where the normalising constant $G^{-1}$ is determined by the
requirement that $\sum_{\vec{n}\in\NN}\pi(\vec{n})=1$.
%  so that
% \begin{equation}
%   \label{eq:xxx}
%   G = \sum_{\vec{n}\in\NN} \prod_{r\in\RR}\frac{\kappa_r^{n_r}}{n_r!}.
% \end{equation}
The simple \emph{product form} of the stationary
distribution~\eqref{eq:4} is a consequence of the fact that the
equations~\eqref{eq:3} \emph{do} have a solution, that is, it is a
consequence of the reversibility of the stationary version of the
process~$\vec{n}(\cdot)$.  Note also that here the stationary
distribution~$\pi$ depends on the parameters~$\nu_r$ and $\mu_r$ only
through their ratios $\kappa_r=\nu_r/\mu_r$, $r\in\RR$.  This result
is not in general true for networks with controls.

In the variation of our model in which calls of each type~$r$ have
holding times which are no longer necessarily exponential
(but with unchanged mean~$\mu_r^{-1}$), it is well-known that the
stationary distribution~$\pi$ of the process $\vec{n}(\cdot)$ continues to
satisfy the detailed balance equations~\eqref{eq:3} and hence also
\eqref{eq:4}.  For a
proof of this \emph{insensitivity} property, see Burman \textit{et
  al.} (1984).

The stationary probability that a call of type~$r$ is \emph{accepted},
is given by
\begin{equation}
  \label{eq:5}
  P_r = \sum_{\vec{n}\in\NN_r} \pi(\vec{n}),
\end{equation}
where $\NN_r=\{\vec{n}\in\NN\colon \vec{n}+\vec{\delta}_r\in\NN\}$.
In Section~\ref{sec:kdr-recursion} we give a recursion which permits a
reasonably efficient calculation of the probabilities~$P_r$ in
networks of small to moderate size.  However, the exact calculation of
acceptance probabilities is usually difficult or impossible in large
networks.  We shall therefore also discuss various approximations.

\subsection{The single resource case}
\label{sec:single-resource}

Consider first the case $\RR=\{1\}$ of a single call type.  For
convenience we drop unnecessary subscripts denoting dependence on
$r\in\RR$; in particular we write $\kappa=\nu/\mu$.  We then have
$\NN=\{n\colon n\le C\}$ for some positive integer~$C$.  The
stationary distribution~$\pi$ is a truncated Poisson distribution, and
the stationary acceptance probability~$P$ is given by Erlang's well-known
formula, that is, by $P=1-\pi(C)=1-E(\kappa,\,C)$, where
\begin{equation}
  \label{eq:6}
  E(\kappa,\,C) = \frac{\kappa^C/C!}{\sum_{n=0}^C\kappa^n/n!}.
\end{equation}
Note also that, from \eqref{eq:2}, the expected number of calls in
progress under the stationary distribution~$\pi$ is given by
$\kappa{}P$.

While exact calculation of blocking probabilities via Erlang's
formula~\eqref{eq:6} is straightforward, it nevertheless provides
insight to give approximations for networks in which $C$ and $\kappa$
are both large.  Formally, we consider the Kelly limiting scheme in
which $C$ and $\kappa$ are allowed to tend to infinity in proportion
to a scale parameter~$N$ with $p=C/\kappa$ held fixed.  The
cases $p>1$, $p=1$ and $p<1$ correspond to the network being, in an
obvious sense, underloaded, critically loaded, and overloaded
respectively.  A relatively straightforward analysis of \eqref{eq:6}
shows that,
\begin{equation}
  \label{eq:7}
  P \to \min(1,p) \qquad\text{as $N\to\infty$.}
\end{equation}
For $p\ge1$ the error in the approximation $P\approx1$ may be estimated by
replacing the truncated Poisson distribution of $n$ by a truncated
normal distribution: for $p>1$ it may be shown to decay at least
exponentially fast in $N$, while for the critically loaded case $p=1$
it may be shown to be $O(N^{-1/2})$ as $N\to\infty$.  For the
overloaded case $p<1$ the approximation $P\approx{}p$ may be refined as
follows.  Observe that in this case, and since $\kappa$ and $C$ are
large, it follows from either \eqref{eq:3} or \eqref{eq:4} that the
stationary distribution of \emph{free capacity} in the network is
approximately geometric 
% with, for $n$ close to $C$,
% \begin{equation}
%   \label{eq:xxx}
%   \pi(C-n) \approx (1-p)p^{n}, \qquad n = 0,1,\dots.
% \end{equation}
and so the stationary expected free capacity is given by the
approximation
\begin{equation}
  \label{eq:8}
  C - \E_\pi(n) \approx \frac{p}{1-p}.
\end{equation}
Combining this with \eqref{eq:2} leads to the very much more refined
approximation for the stationary acceptance probability given by
\begin{equation}
  \label{eq:9}
  P \approx p - \frac{p}{\kappa(1-p)}.
\end{equation}
The error in this approximation may be shown to be $o(N^{-1})$ as
$N\to\infty$, so also that in the original approximation $P\approx{}p$
is $O(N^{-1})$.

\subsection{The Kaufman-Dziong-Roberts (KDR) recursion}
\label{sec:kdr-recursion}

For the general model of an uncontrolled network, we now take the set
$\NN$ to be given by a set of capacity constraints of the
form~\eqref{eq:1}.  We give here an efficient recursion for the
determination of stationary acceptance probabilities, due in the case
$\JJ=\{1\}$ to Kaufman (1981) and in the general case to
Dziong and Roberts (1987).

% For each time~$t$, define $\vec{m}(t)=(m_j(t),\,j\in\JJ)$ where
% \begin{equation}
%   \label{eq:39}
%   m_j(t) = \sum_{r\in\RR}A_{jr}n_r(t).
% \end{equation}
% The process $\vec{m}(\cdot)$ thus describes the resource usage in the
% network.  It takes values in the set
% $\MM=\{\vec{m}\in\Zp^J\colon0\le{}m_j\le{}C_j,\;j\in\JJ)$ where we
% write $J=|\JJ|$.  The stationary distribution~$\pi'$ of this process
% is given by
% \begin{equation}
%   \label{eq:10}
%   \pi'(\vec{m}) = \sum_{\vec{n}\colon A\vec{n}=\vec{m}} \pi(\vec{n}),
%   \qquad \vec{m}\in\MM.
% \end{equation}

Recall that $\pi'$ is the stationary distribution of the
process~$\vec{m}(\cdot)$ defined in the Introduction.  Since a call of
type~$r$ arriving at time~$t$ is accepted if and only if
$m_j(t-)+A_{jr}\le{}C_j$ for all $j$ such that $A_{jr}\ge1$ (where
$\vec{m}(t-)$ denotes the state of the process~$\vec{m}(\cdot)$
immediately prior to the arrival of the call), it follows
that a knowledge of $\pi'$ is sufficient to determine stationary
acceptance probabilities.  Typically the size~$J$ of the set~$\JJ$ is
smaller than the size~$R$ of the set~$\RR$, and so the dimension of
the space~$\MM$ defined by~\eqref{eq:46} is smaller than that of
$\NN$.  Thus a direct calculation of $\pi'$, avoiding that of $\pi$,
is usually much more efficient for determining acceptance
probabilities.

For each $r\in\RR$, define the vector $\vec{A}_r=(A_{jr},\,j\in\JJ)$.
For each $\vec{m}\in\MM$ and $r\in\RR$, summing the detailed balance
equations~\eqref{eq:3} over $\vec{n}$ such that $A\vec{n}=\vec{m}$
and using also \eqref{eq:44} yields
\begin{equation}
  \label{eq:11}
  \kappa_r \pi'(\vec{m}-\vec{A}_r) = \E(n_r \gvn \vec{m})\pi'(\vec{m}),
  \qquad r\in\RR, \quad \vec{m}\in\MM,
\end{equation}
where
\begin{displaymath}
  \E(n_r \gvn \vec{m})
  = \frac{\sum_{\vec{n}\colon A\vec{n}=\vec{m}} n_r\pi(\vec{n})}%
  {\sum_{\vec{n}\colon A\vec{n}=\vec{m}} \pi(\vec{n})}
\end{displaymath}
is the stationary expected value of $n_r$ given $A\vec{n}=\vec{m}$.
Since, for each $\vec{m}$ and each $j$, we have
$\sum_{r\in\RR}A_{jr}\E(n_r \gvn \vec{m})=m_j$, it follows
from~\eqref{eq:11} that
\begin{equation}
  \label{eq:12}
  \sum_{r\in\RR}A_{jr}\kappa_r \pi'(\vec{m}-\vec{A}_r)
  = m_j\pi'(\vec{m}),
  \qquad \vec{m}\in\MM, \quad j\in\JJ.
\end{equation}
This is the Kaufman-Dziong-Roberts (KDR) recursion on the set $\MM$,
enabling the direct determination of successive values of
$\pi'(\vec{m})$ as multiples of $\pi'(\vec{0})$.  The entire
distribution $\pi'$ is then determined uniquely by the requirement
that $\sum_{\vec{m}\in\MM}\pi'(\vec{m})=1$.

\subsection{Approximations for large networks}
\label{sec:appr-large-netw}

We now suppose that $\kappa_r$, $r\in\RR$, and  $C_j$, $j\in\JJ$, are
sufficiently large that the exact calculation of the stationary
distributions~$\pi$ or $\pi'$ is impracticable.  We seek good
approximations for the latter and for acceptance probabilities.

\paragraph{A simple approximation}
\label{sec:simple-approximation}

We give first a simple approximation, due to Kelly (1986), which
generalises the approximation $P\approx\min(1,p)$ of
Section~\ref{sec:single-resource} for the single-resource case.  To
provide asymptotic justification we again consider the Kelly limiting
scheme, in which the parameters~$\kappa_r$ and $C_j$ are allowed to
increase in proportion to a scale parameter~$N$, the sets $\RR$, $\JJ$ and
the matrix~$A$ being held fixed.  We assume (this is largely for
simplicity) that the matrix~$A$ is such that for each $\vec{m}\in\MM$
there is at least one $\vec{n}\in\Z$ such that $A\vec{n}=\vec{m}$.
(This implies in particular that the matrix~$A$ is of full rank.)  We
outline an argument based
on the equations~\eqref{eq:11} and the KDR recursion~\eqref{eq:12}.

Suppose that $\pi'(\vec{m})$ is maximised at $\vec{m}^*\in\MM$.  The
distribution~\eqref{eq:4} of $\pi$ is a truncation of a product of
independent Poisson distributions each of which has a standard
deviation which is $O(N^{1/2})$ as the scale parameter~$N$ increases.
{}From this and from the mapping of $\pi$ to $\pi'$, it follows that
all but an arbitrarily small fraction of the distribution of $\pi'$ is
concentrated within a region~$\MM^*\subseteq\MM$ such that the
components of $\vec{m}\in\MM^*$ differ from those of $\vec{m}^*$ by an
amount which is again $O(N^{1/2})$ as $N$ increases.  Further, it is
not too difficult to show from the above condition on the matrix~$A$
that, for each $r$, $\E(n_r\gvn\vec{m})$ varies smoothly with
$\vec{m}$, and that within $\MM^*$ we may make the approximation
$\E(n_r\gvn\vec{m})\approx\E(n_r\gvn\vec{m}^*)$ (the error yet again
being $O(N^{1/2})$ as $N$ increases).  It now follows from
\eqref{eq:11} that within $\MM^*$ we have
\begin{equation}
  \label{eq:13}
  \pi'(\vec{m}) \approx \pi'(\vec{m}^*)\prod_{j\in\JJ}p_j^{m^*_j-m_j},
\end{equation}
where
% , for each $r$, we have
% \begin{displaymath}
%   \frac{\E(n_r\gvn\vec{m}^*)}{\kappa_r} = \prod_{j\in\JJ}p_j^{A_{jr}},
% \end{displaymath}
necessarily, since $\vec{m}^*$ maximises $\pi'(\vec{m})$,
\begin{align}
  0 \le p_j & \le 1, \qquad j\in\JJ  \label{eq:14},\\
  p_j & = 1, \qquad\text{for $j$ such that $m^*_j<C_j$  \label{eq:15}.}
\end{align}
Further, from \eqref{eq:12},
\begin{equation}
  \label{eq:16}
  \sum_{r\in\RR}A_{jr}\kappa_r\prod_{k\in\JJ}p_k^{A_{kr}} = m^*_j \le C_j,
  \qquad j\in\JJ.
\end{equation}

Thus, from \eqref{eq:13}, within
$\MM^*$ the stationary distribution~$\pi'$ of $\vec{m}(\cdot)$
does indeed have the approximate factorisation~\eqref{eq:45}, where
each of the component distributions~$\pi'_j$ is here geometric (and
where in the case $p_j=1$ the geometric distribution becomes uniform).
Further, for each $r$ and each $j$, we have
\begin{displaymath}
  \pi'_j(\{m_j\colon m_j\le{}C_j-A_{jr}\})\approx{}p_j^{A_{jr}}.
\end{displaymath} 
Thus the stationary acceptance probabilities $P_r$ are given by the
approximation
\begin{equation}
  \label{eq:18}
  P_r \approx \prod_{j\in\JJ}p_j^{A_{jr}}, \qquad r\in\RR.
\end{equation}
Kelly (1986) considered an optimisation problem from which it
follows that the equations~\eqref{eq:14}--\eqref{eq:16} determine the
vectors~$\vec{m}^*$ and $\vec{p}=(p_j,\,j\in\JJ)$ uniquely.
% of the optimisation problem
% \begin{displaymath}
%   \text{minimise\ } f(\vec{p})
%   := \sum_{r\in\RR}\kappa_r\prod_{j\in\JJ}p_j^{A_{jr}}
%   - \sum_{j\in\JJ}C_j\log p_j,
%   \qquad\text{subject to 
%   $p_j\in(0,1]$, ~ $j\in\JJ$.}
% \end{displaymath}
% Since the objective function $f$ is convex, it follows that
% \eqref{eq:14}--\eqref{eq:16} determine $\vec{p}$ uniquely.  Kelly
He further showed, in an approach based on consideration of the
stationary distribution~$\pi$, that the approximation~\eqref{eq:18}
becomes exact as the scale parameter~$N$ tends to infinity.

\paragraph{A refined approximation}

The \emph{(multiservice) reduced load} or \emph{knapsack
  approximation} (Dziong and Roberts, 1987, see also
Ross, 1995) is a more refined approximation than that defined
above.  It is given by retaining the approximate
factorisation~\eqref{eq:45} of the stationary distribution~$\pi'$ of
$\vec{m}(\cdot)$.  However, subject to this assumed factorisation, the
estimation of the component distributions~$\pi'_j$ is refined.

For each~$j\in\JJ$ and $r\in\RR$, define
\begin{equation}
  \label{eq:19}
  p_{jr} = \sum_{m_j=0}^{C_j-A_{jr}} \pi'_j(m_j);
\end{equation}
note that $p_{jr}=1$ if $A_{jr}=0$.  For fixed $j$, substitution of
\eqref{eq:45} into the KDR recursion~\eqref{eq:12} and summation over
all $m_k$ for all $k\neq{}j$ yields
\begin{equation}
  \label{eq:20}
  \sum_{r\in\RR}A_{jr}\biggl(\kappa_r \prod_{k\neq{}j}p_{kr}\biggr)%
  \pi'_j(m_j-A_{jr})
  = m_j\pi'_j(m_j),
  \qquad 1 \le m_j \le C_j, \quad j\in\JJ
\end{equation}
(where, as usual, we make the convention $\pi'_j(m_j)=0$ for $m_j<0$).
This is the one-dimensional KDR recursion associated with a single
resource constraint~$j$, and is readily solved to determine $\pi'_j$
and hence the probabilities~$p_{jr}$, $r\in\RR$, in terms of the
probabilities~$p_{kr}$, $r\in\RR$, for all $k\neq{}j$.  We are thus
led to a set of fixed point equations in the probabilities~$p_{jr}$,
for which the existence---but not always the uniqueness, see Chung and
Ross (1993)---of a solution is guaranteed.  From~\eqref{eq:45},
the probability that a call of type~$r$ is accepted is then given by
\begin{equation}
  \label{eq:21}
  P_r = \prod_{j\in\JJ}p_{jr}.
\end{equation}
We remark that the recursion~\eqref{eq:20} corresponds to a modified
network in which there is a single resource constraint~$j$ and each
arrival rate~$\kappa_r$ is reduced to $\kappa_r
\prod_{k\neq{}j}p_{kr}$.  This \emph{reduced load approximation} is of
course exact in the case of a single-resource network.

In the case where each $A_{jr}$ can only take the values $0$ or $1$ we
may set $p_j=p_{jr}$ for $r$ such that $A_{jr}=1$.  The fixed point
equations~\eqref{eq:19} and \eqref{eq:20} then reduce to
\begin{equation}
  \label{eq:22}
  p_j = 1 - E\Biggl(\sum_{r\in\RR}\kappa_r\prod_{k\neq{}j}p_k^{A_{kr}},\,C_j\Biggr)
\end{equation}
where $E$ is the Erlang function~\eqref{eq:6}.  This case is
the well-known \emph{Erlang fixed point approximation} (EFPA) and has
a unique solution, see Kelly (1986), and also
Ross (1995).  It yields acceptance probabilities which are
known to be asymptotically exact in the Kelly limiting scheme
discussed above, and also, under appropriate conditions, in the
\emph{diverse routing} limit discussed in the Introduction---see Whitt (1985), 
and Ziedins and Kelly (1989).  The EFPA also has
an extension to the case of general $A_{jr}$, which may be regarded as
a simplified version of the reduced load approximation.  As with the
latter approximation the EFPA may here have multiple solutions.

\section{Controlled loss networks: stationary behaviour}
\label{sec:contr-loss-netw}

We now study the more general version of a loss network, in
which calls are subject to acceptance controls, and the issues are
those of achieving optimal performance.

\subsection{Single resource networks}
\label{sec:simple-model}

We consider a simple model which illustrates some ideas of optimal
control---in particular those of \emph{robustness} of the control
strategy with respect to variations in arrival rates (which may in
practice be unknown, or vary over time).

Suppose that $\RR=\{1,2\}$ and that as usual calls of each type~$r$
arrive at rate $\nu_r$ and have holding times which are exponentially
distributed with mean $\mu_r^{-1}$.  Suppose further that there is a
single resource of capacity~$C$ and that a call of either type requires 
one unit of this capacity, so that the constraints~\eqref{eq:1} here 
reduce to $n_1 + n_2 \le C$.
We assume that calls of type~$1$ have greater value per
unit time than those of type~$2$, so that it is desirable to choose
the acceptance regions~$\AA_r$, $r=1,2$, so as to maximise the linear
function
\begin{equation}
  \label{eq:23}
  \phi(P_1,P_2) := a_1\kappa_1P_1 + a_2\kappa_2P_2,
\end{equation}
for some $a_1>a_2>0$ (where, again as usual, for each call type~$r$,
$\kappa_r=\nu_r/\mu_r$ and $P_r$ is the stationary
acceptance probability.)  An upper bound
for the expression in \eqref{eq:23} is given by the solution of the
linear programming problem, in the variables~$P_1$, $P_2$,
\begin{equation}
  \label{eq:40}
  \text{maximise $\phi(P_1,P_2)$,\quad subject to $P_r\in[0,1]$ for
    $r=1,2$, \quad $\kappa_1P_1+\kappa_2P_2\le{}C$} 
\end{equation}
(where the latter constraint follows from~\eqref{eq:2}).  It is easy
to see that the solution of this problem is characterised uniquely by
the conditions
\begin{alignat}{2}
  P_1 & = P_2 = 1, && \quad\text{whenever $\kappa_1P_1+\kappa_2P_2 < C$,}
  \label{eq:41}\\
  P_2 & = 0, && \quad\text{whenever $P_1 < 1$.} \label{eq:42}
\end{alignat}
% An elementary calculation therefore shows that this upper bound is
% given by taking
% \begin{equation}
%   \label{eq:24}
%   P_1 = \min\left(1,\frac{C}{\kappa_1}\right),
%   \qquad
%   P_2 = \min\left(1,\frac{C-\kappa_1P_1}{\kappa_2}\right).
% \end{equation}
It is clearly not possible to choose the acceptance regions~$\AA_1$,
$\AA_2$ so that the corresponding values of $P_1$, $P_2$ solve exactly
the problem~\eqref{eq:40}.  However, we show below
that this solution may be achieved asymptotically as the size of the
system is allowed to increase, and further that there is an 
asymptotically optimal control that is
both simple and robust with respect to variations in the parameters
$\kappa_1$, $\kappa_2$.

We consider first the form of the optimal control in the special case
$\mu_1=\mu_2$.  Here since, at the arrival time~$t$ of any call, those
calls already within the system are indistinguishable with respect to
type, it is clear that the optimal decision on call admission is a
function only of the arriving call type and of the total volume
$m(t-)=n_1(t-)+n_2(t-)$ of calls already in the system.  A formal proof
is a straightforward exercise in Markov decision theory.  Further,
simple coupling arguments show that, for an incoming call of either
type arriving at time~$t$ and any $0<m<C$, if it is advantageous to
accept the call when $m(t-)=m$, then it is also advantageous to accept
the call when $m(t-)=m-1$.  It follows that the optimal acceptance
regions are of the form
\begin{align}
  \AA_1 & = \{\vec{n}\colon n_1+n_2 < C\} \label{eq:25}\\
  \AA_2 & = \{\vec{n}\colon n_1+n_2 < C-k\} \label{eq:26}
\end{align}
for some \emph{reservation parameter}~$k$, whose  optimal value
depends on $C$, $\kappa_1$ and $\kappa_2$.

Consider now the general case where we do not necessarily have
$\mu_1=\mu_2$, and suppose that $C$, $\nu_1$ and $\nu_2$ are
large.  More formally we again have in mind the Kelly limiting scheme
in which these parameters are allowed to increase in proportion to some
scale parameter~$N$ which tends to infinity (while $\mu_1$, $\mu_2$
are held fixed).  We further suppose that
the acceptance regions are again as given by~\eqref{eq:25} and
\eqref{eq:26}, where the reservation parameter~$k$ increases slowly
with $N$, i.e.\ in such a way that
\begin{equation}
  \label{eq:27}
  k \to \infty, \qquad k/C \to 0, \qquad \text{as $N\to\infty$}.
\end{equation}
It is convenient to let $P_1$, $P_2$ denote the limiting acceptance
probabilities.  In the case $\kappa_1+\kappa_2\le{}C$, it is not
difficult to see that, since $k/C\to0$ as $N\to\infty$, we have
$P_1=P_2=1$, so that $P_1$, $P_2$ solve the optimisation
problem~\eqref{eq:40}.  Consider now the case $\kappa_1+\kappa_2>C$.
Here, again since $k/C\to0$ as $N\to\infty$, it follows that, in the
limit, the capacity of the network is fully utilised.  Further, if
$\kappa_1$ is sufficiently large that $P_1<1$ (informally, even for
large $N$, calls of type~$1$ are being rejected in significant
numbers), then the effect of the increasing reservation parameter~$k$ is
such that, again in the limit, the network remains sufficiently close
to capacity to ensure that no calls of type~$2$ are accepted, and
hence $P_2=0$.  It now follows that when $\kappa_1+\kappa_2>C$, the
limiting acceptance probabilities $P_1$, $P_2$ satisfy the
conditions~\eqref{eq:41} and \eqref{eq:42} and so again solve the
optimisation problem~\eqref{eq:40}.

The above analysis demonstrates the asymptotic optimality of any
strategy based on the use of a reservation parameter~$k$, provided
only that, in the limiting regime, $k$ increases in accordance with
\eqref{eq:27}.  In practice, in a large network (here for large~$C$),
only a small value of $k$ is required in order to achieve optimal
performance.  We also observe that the performance of a reservation
parameter strategy is indeed robust with respect to variations in
$\kappa_1$, $\kappa_2$.

This analysis also extends easily to the case where there are more
than two call types, and also, with a little more difficulty, to that
where the capacity constraint is of the form
$\sum_{r\in\RR}A_r{}n_r\le{}C$ for general positive integers~$A_r$
(see Bean \textit{et al.}, 1995).
Here a different reservation parameter may be used for each call type,
and, in the Kelly limiting scheme, a complete prioritisation and
optimal control are again achieved asymptotically by allowing the
differences between the reservation parameters to increase slowly.

% \sz{To what extent does this result (asymptotic achievement in the
%   Kelly limit of the upper bound given by the simple optimisation
%   problem, by the use of an appropriate reservation-based control)
%   extend to the case of multiple resource constraints.  This surely
%   should be in the literature---but where?}

\subsection{Multiple resource models}
\label{sec:mult-reso-models}

Consider now the general case of the canonical model in which there is a
set of resources~$\JJ$ and in which state $\vec{n}$ of the network is
subject to the constraints~\eqref{eq:1}.  Suppose that it is again
desirable to choose admission controls so as to maximise the
linear function $\phi(\vec{P}):=\sum_{r=1}^Ra_r\kappa_rP_r$ of the
stationary acceptance probabilities~$P_r$, for given constants~$a_r$,
$r\in\RR$.   As in Section~\ref{sec:simple-model}, we may consider the
linear programming problem
\begin{equation}
  \label{eq:37}
  \text{maximise $\phi(\vec{P})$,
    \quad subject to $P_r\in[0,1]$ for
    $r\in\RR$,
    \quad $\sum_{r=1}^RA_{jr}\kappa_rP_r\le{}C_j$ for $j\in\JJ$,} 
\end{equation}
which provides an upper bound on the achievable values of the
objective function~$\phi$.  It is easy to see that this value may be
asymptotically achieved within the Kelly limiting regime by reserving
capacity~$A_{jr}\kappa_rP_r$ at each resource~$j$ solely for calls of
each type~$r$, where here $\vec{P}$ is the solution of the
problem~(\ref{eq:37}).  However this strategy is neither optimal in
networks of finite capacity, nor is it robust with respect to
variations in the parameters~$\kappa_r$.  At the opposite end of
the spectrum from this complete partitioning policy is that of complete
sharing.  The latter can lead to unfairness if there are asymmetric
traffic patterns, with the potential for some call types to receive
better service than others.  In practice it is expected
that good strategies will be based on the sharing of resources and the
use of reservation parameters---as was shown to be optimal for single
resource networks in Section~\ref{sec:simple-model}.  

In the case of communications networks it is natural to allow also
\emph{alternative routing}, as described in the Introduction.  An
upper bound for the achievable performance is given by supposing that
\emph{repacking} is possible, i.e.\ that calls in progress may be
rerouted as necessary.  In this case, our model for the network 
reduces to an instance of the canonical model (as defined in the
Introduction)
with appropriately redefined set~$\JJ$, matrix~$A=(A_{jr})$ and
capacities~$C_j$.  The upper bound on $\phi(\vec{P})$ given by
the linear programming problem~\eqref{eq:37} is then also an upper bound in
the more usual case in which repacking is not allowed.  In the latter
case practical control strategies are again based on the use of appropriate
reservation parameters, and there is some hope that performance close to
the upper bound above may be achieved in networks with
sufficiently large capacities or sufficient diversity of routing,
even without repacking.  In applications reservation parameters are
generally used
to prioritise different traffic streams.  In networks with
alternative routing they also prevent the occurrence of network instabilities, 
where, for fixed parameter values, the
network may have two or more relatively stable operating regimes---one
in which most calls are directly routed, and others in which many
calls are alternatively routed, with a resulting severe
degradation of performance (see Gibbens {\it et al.}, 1990, Kelly, 1991).
By giving priority to directly routed traffic, the use of reservation
parameters prevents the network
from slipping into an inefficient operating state.

There have been numerous investigations of control strategies for
communications networks
that employ either fixed or, particularly, alternative routing. 
% with the use of reservation parameters to prevent instability and
% achieve an efficient use of resources.  
Such strategies are
often studied in the context of fully connected networks.  
Two of the most commonly studied are 
{\em least busy alternative} (LBA) routing
and \emph{dynamic alternative routing} (DAR).
LBA routing seeks to route calls directly if possible, and otherwise
routes them via that path which minimises the maximum occupancy on
any of its links.  Directly routed calls are usually ``protected'' with
some form of reservation parameter (Kelly, 1991, Marbukh, 1993).  Hunt and
Laws (1993) showed that, for fully connected networks which permit only
two-link alternative routes, LBA routing is asymptotically optimal 
in the diverse routing limit  (see
Section~\ref{sec:diverse-rout-limit}).  
This policy is
robust to changes in traffic patterns, but has the difficulty that
it requires information on the current states of all possible alternative
paths before an alternative routing decision is made.  

%One way of
%overcoming this is to take an approximate version of LBA, such as the
%aggregated least busy alternative (ALBA) routing scheme outlined
%in Mitra, Gibbens and Huang (1993) and Mitra and Gibbens(1992), which
%maintains information about the occupancy of a resource in a much
%reduced fashion (rather than the exact occupancy, only the approximate level
%is recorded, where there may be as few as two or three possible levels).
A much simpler routing scheme is DAR (Gibbens {\it et al.}, 1989,
Gibbens and Kelly, 1990).
In this scheme, for each pair of nodes, a record is maintained
of the current preferred alternative route, and this is the
one that is used if a call cannot be routed directly.  If neither
the direct route nor the current preferred alternative route are
available, then the call is rejected, and a new preferred alternative
route is chosen at random from those available.  Directly routed
traffic is again usually protected by a reservation parameter.
This policy is easy to implement.  It does not require 
information about the current state of the system to be held at any
node, just a
record of the current preferred alternative route to other nodes.
It is also robust to changes in traffic patterns---alternative routes
on which the load increases will be discarded and replaced by routes
on which the load is lower.  Neither LBA routing nor DAR require traffic
rates to be known or estimated (except approximately, in order to set
the appropriate level of the reservation parameters).

Acceptance probabilities for controlled loss networks are usually
estimated using a generalised version of the reduced load or knapsack
approximation of Section~\ref{sec:appr-large-netw}.  As there, we make
the approximation~\eqref{eq:45} for the stationary distribution~$\pi'$
of the resource occupancy process~$\vec{m}(\cdot)$.  Each of the
marginal distributions~$\pi'_j$ is estimated as the stationary
distribution of a Markov process on $\{0,\dots,C_j\}$ which
approximates the behaviour of the resource~$j$ considered in
isolation.  Let $p_{jr}$ be the probability under this distribution
that a call of type~$r$ is accepted, subject to the controls of the
model, with $p_{jr}=1$ if $A_{jr}=0$.  In the case of the canonical
model, in which no alternative resource usage is allowed, calls of
each type~$r$ are assumed to arrive at resource~$j$ at a rate
$\nu_r\prod_{k\neq{}j}p_{kr}$---this is the ``reduced load'' for calls
of type~$r$ at this resource; further, calls of this type arriving at
this resource are subject to the acceptance controls of the model and,
if accepted, depart at rate~$\mu_r$ as usual.  The estimated
stationary distribution~$\pi'_j$ then determines the acceptance
probabilities~$p_{jr}$ at the resource~$j$.  Thus we are again led to
a set of fixed point equations which determine---not always
uniquely---the acceptance probabilities $p_{jr}$ for all $r\in\RR$ and
$j\in\JJ$.  Finally the stationary network acceptance
probability~$P_r$ for calls of each type~$r$ is again given by
$P_r=\prod_{j\in\JJ}p_{jr}$.

In the case of a communications network where the canonical model is
extended by allowing the possibility of alternative routing, it is
necessary to modify the above approximation.  Suppose, for example,
that a link (resource)~$j$ forms part of the second choice route for
calls of type~$r$.  Then, in the one-dimensional process associated
with link~$j$, the arrival rate for calls of type~$r$ is taken to be
the product of the arrival rate~$\nu_r$ at the network, the
probability that a call of this type is rejected on its first-choice
route, and (as before) the probabilities that the call can be accepted
at each of the remaining resources on the alternative route (see e.g.\
Gibbens and Kelly, 1990).

The basis of the reduced load approximation is the approximate
factorisation of the distribution~$\pi'$ above.  In the case of
controlled networks, this approximation fails to become exact under
the Kelly limiting regime in which capacities and arrival rates increase
in proportion.  It may, however, be expected to hold under
sufficiently diverse routing.  It is known to be remarkably accurate
in most applications.

\section{Dynamical behaviour and stability}
\label{sec:dyn}

\subsection{Fluid limits for large capacity networks}
\label{sec:fluid-limits-large}

We now consider the dynamical behaviour of large networks.  As well as
such behaviour being of interest in its own right---for example in
networks in which input rates change suddenly, \emph{fixed points} of
network dynamics correspond to equilibrium, or quasi-equilibrium,
states of the network (see below).  The identification of such points
is often the key to understanding long-term behaviour, in particular
to resolving stability questions and determining
stationary distributions where (as is usual) the latter may not be
directly calculated.  However, we note that it is characteristic of
loss networks that, from any initial state, equilibrium is effectively
achieved within a very few call holding times, so that transient
\emph{performance} is of less significance than is the case for
networks which permit queueing.

We describe a theory first suggested by Kelly (1991).  We yet again assume the
Kelly limiting scheme described in the Introduction, in which the
network topology is held fixed and arrival rates and capacities are
allowed to increase in proportion.  More explicitly, we consider a
sequence of networks satisfying our usual Markov assumptions (though
this is not strictly necessary) and indexed by a scale parameter~$N$.
All members of the sequence are identical in respect of the (finite)
sets~$\RR$, $\JJ$, the matrix $A=(A_{jr},j\in\JJ,\,r\in\RR)$, and the
departure rates $\mu_r$, $r\in\RR$.  For the $N$th member of the
sequence, calls of each type~$r$ arrive at rate $N\nu_r$ for some
vector of parameters~$\vec{\nu}$, and the capacity of each
resource~$j$ is $NC_j$ for some vector of parameters~$\vec{C}$,
where, for simplicity, we take each $C_j$ to be integer-valued.  As
always, it is convenient to define $\kappa_r=\nu_r/\mu_r$ for each
$r\in\RR$.

We now describe the rules whereby calls are accepted.  For each $N$,
let $\vec{n}^N(t)=(n^N_r(t),~r\in\RR)$, where $n^N_r(t)$ is the number
of calls of type~$r$ in progress at time~$t$.  Define also the
\emph{free capacity} process
$\fc{\vec{m}}^N(\cdot)=(\fc{m}^N_j(\cdot),~j\in\JJ)$ where each
$\fc{m}^N_j(t)=NC_j-\sum_{r\in\RR}A_{jr}n^N_r(t)$ is the free capacity
of resource~$j$ at time~$t$.
% We formally regard the
% process~$\fc{m}^N(\cdot)$ as taking values in the compactified
% space~$E=(\Zp\cup\{\infty\})^J$, where $J=|{\cal J}|$.  
A call of type~$r$ arriving at time $t$ is accepted if and only if the
free capacity $\fc{\vec{m}}^N(t-)$ of the system, immediately prior to
its arrival, belongs to some acceptance
region~${\fc{\AA}}_r\subset\Zp^J$.
% whose indicator function is required to satisfy the somewhat technical
% requirement of continuity with respect to the topology on $E$ given by
% the product of the one-point compactifications of $\Zp$ (a condition
% which will always be satisfied in practice).  
We take the acceptance regions~${\fc{\AA}}_r$, $r\in\RR$, to be
independent of $N$, although, in a refinement of the theory, some
dependence may be allowed.  Note that, in a change from our earlier
conventions, the acceptance regions~${\fc{\AA}}_r$ are defined in
terms of the \emph{free} capacity of each system.

While the above description defines instances of the canonical model
of the Introduction, more sophisticated controls, such as those
involving the use of alternative routing in communications networks,
may be modelled by the suitable redefinition of input streams and
acceptance sets (see Hunt and Kurtz, 1994).

For each $N$, define the normalised
process~$\vec{x}^N(\cdot)=\vec{n}^N(\cdot)/N$, which takes values in
the space
\begin{equation}
  \label{eq:30}
  X=\{\vec{x}\in\Rp^R\colon
  \text{$\sum_{r\in\RR}A_{jr}x_r\le C_j$ for all~$j\in\JJ$}\}.
\end{equation}
Assume that, as $N\to\infty$, the initial state~$\vec{x}^N(0)$
converges in distribution to some $\vec{x}(0)\in{}X$, which, for
simplicity, we take to be deterministic.  Then we might expect that
the process~$\vec{x}^N(\cdot)$ should similarly converge in
distribution to a \emph{fluid limit} process~$\vec{x}(\cdot)$ taking
values in the space $X$,
% We show informally
% that this is the case, where the dynamics of $x(\cdot)$ are 
with dynamics given by
\begin{equation}
  \label{eq:31}
  x_r(t) = x_r(0) + \int_0^t (\nu_r \tilde{P}_r(u) - \mu_r x_r(u))du,
  \qquad r \in \RR,
\end{equation} 
where, for each $t$, $\tilde{P}_r(t)$ corresponds to the limiting rate at which
calls of each type $r$ are being accepted at time~$t$.  

A rigorous convergence result is given by Hunt and Kurtz (1994).  A
somewhat technical condition (always likely to be satisfied in
applications) is required on the acceptance sets~$\fc{\AA}_r$.
However, the main difficulty is that in some, usually rather
pathological, cases the limiting acceptance rates $\tilde{P}_r(t)$ may fail to
be unique.

In many cases, though, it is possible to show that, for each $r$,
there does exist a unique function $P_r$ on $X$ such that, for each
$t$, we have $\tilde{P}_r(t)=P_r(\vec{x}(t))$.
% We may then
% define a \emph{velocity field}~$\vec{v}=(v_r,~r\in\RR)$ on $X$ by
% \begin{equation}
%   \label{eq:xxx}
%   v_r(\vec{x}) = \nu_rP_r(\vec{x}) - \mu_r x_r,
%   \qquad r\in\RR,
% \end{equation}
% so that \eqref{eq:31} then becomes
% \begin{equation} 
%   \label{eq:xxx}
%   x_r(t) = x_r(0) + \int_0^t v_r(\vec{x}(u))du.
%   \qquad r \in \RR,
% \end{equation} 
In general, the trajectories of the limit process~$\vec{x}(\cdot)$ are
then deterministic functions of their initial positions $\vec{x}(0)$.
The fixed points $\hat{\vec{x}}$ of the limit process $\vec{x}(\cdot)$ are
given by the solutions of
\begin{equation}
  \label{eq:32}
  \nu_r P_r(\hat{x}) = \mu_r \hat{x}_r, \qquad r\in\RR.
\end{equation}
In the case of a single fixed point $\hat{\vec{x}}$, to which all
trajectories of $\vec{x}(\cdot)$ converge, it may be shown that the
stationary distribution of the original normalised
process~$\vec{x}^N(\cdot)$ converges to that concentrated on the
single point~$\hat{\vec{x}}$.  Then in particular, for each~$r$,
$P_r(\hat{\vec{x}})$ is the limiting stationary acceptance probability
for calls of type~$r$.  In the case of multiple fixed points, those
which are locally stable correspond to ``quasi-stationary''
distributions of the process~$\vec{x}^N(\cdot)$, i.e.\ regimes which
are maintained over periods of time which are lengthy but finite.

\subsection{Single resource networks}
\label{sec:single-reso-netw}

As the simplest non-trivial application of the above theory, we
consider the case $J=1$ of a single resource, for which equilibrium
behaviour was described in Section~\ref{sec:simple-model}.  It is
again convenient to write $A_r$ for $A_{1r}$ for each $r$, and
similarly $C$ for $C_1$.  The technical condition referred to above on
the acceptance sets~$\fc{\AA}_r\subseteq\Zp$, here reduces to the
requirement that, for each $r$, \emph{either} $m\in\fc{\AA}_r$ for all
sufficiently large $m\in\Zp$---we let $\RR^*$ denote the set of such
$r$---\emph{or} $m\notin\fc{\AA}_r$ for all sufficiently large
$m\in\Zp$.

Here the functions~$P_r$ defined above always exist (see Hunt and
Kurtz, 1994).  To identify them, define, for each $\vec{x}\in{}X$, the
Markov process $\fc{m}_\vec{x}(\cdot)$ on $\Zp$ with transition rates
given by
\begin{equation}
  \label{eq:33}
  \fc{m}\rightarrow
  \begin{cases}
    \fc{m}-A_r & \text{at rate $\nu_rI_{\{\fc{m}\in\fc{\AA}_r\}}$}\\
    \fc{m}+A_r & \text{at rate $\mu_rx_r$,}
  \end{cases}
\end{equation}
Let $\pi_\vec{x}$ be the stationary distribution of this
process where it exists.  Define $\bar{X}\subseteq{}X$ by
\begin{equation}
  \label{eq:34}
  \bar{X} = \{\vec{x} \in X\colon
  \sum_{r\in\RR}A_rx_r=C \text{ and $\pi_\vec{x}$ exists}\}. 
\end{equation}
(The set $\bar{X}$ may be thought of as consisting of those points in
$X$ for which the limiting dynamics are ``blocking''.)  Then, for
$\vec{x}\in\bar{X}$, we have $P_r(\vec{x})=\pi_\vec{x}(\fc{\AA}_r)$
for all $r$; for $\vec{x}\in{}X\setminus\bar{X}$, we have
$P_r(\vec{x})=1$ for $r\in\RR^*$ and $P_r(\vec{x})=0$ for
$r\notin\RR^*$.  The fixed points~$\hat{\vec{x}}$ of the
limiting dynamics (in general there may be more than one such) are
then given by the solutions of~\eqref{eq:32}.

Consider now the case of reservation-type controls, and suppose that
the call types are arranged in order of decreasing priority.  The
acceptance regions are thus given by
$\fc{\AA}_r=\{\fc{m}\colon\fc{m}\ge{}k_r+A_r\}$ for some
$0=k_1\le{}k_2\le\dots\le{}k_R$ and we have $\RR^*=\RR$.  It is easy
to see that, in the \emph{light traffic} case given by
$\sum_{r\in\RR}A_r\kappa_r\le{}C$, the single fixed point~$\hat{\vec{x}}$ of
the limiting dynamics is given by $\hat{x}_r=\kappa_r$ for all $r$, and that
all trajectories of these dynamics converge to $\hat{\vec{x}}$.
In the \emph{heavy traffic} case given by
$\sum_{r\in\RR}A_r\kappa_r>C$, define $\hat{X}\subseteq{}X$ by
\begin{displaymath}
  \hat{X} = \{\vec{x} \in X\colon \sum_{r\in\RR}A_rx_r=C
  \text{ and $x_r<\kappa_r$ for all $r\in\RR$}\}.
\end{displaymath}
Then it is straightforward to show that $\hat{X}\subseteq\bar{X}$ and
that all fixed points of the limiting dynamics lie within $\hat{X}$
(see Bean \textit{et al.}, 1995).
In the case where $A_r=1$ for all $r$, it is also straightforward to
show that there is a unique fixed point.  
% To see this latter result, suppose that
% $\vec{x}^{(1)},\vec{x}^{(2)}\in\hat{X}$ are both fixed points and
% suppose, without loss of generality, that
% $\sum_{r\in\RR}\mu_rx^{(1)}_r\ge\sum_{r\in\RR}\mu_rx^{(2)}_r$.  Then
% consideration of the processes~$\pi_{\vec{x}^{(1)}}$,
% $\pi_{\vec{x}^{(2)}}$ defined by \eqref{eq:33} shows that
% $P_r(\vec{x}^{(1)})\ge{}P_r(\vec{x}^{(2)})$ for all $r$, and so it
% follows from \eqref{eq:32} that $x^{(1)}_r\ge{}x^{(2)}_r$ for all $r$.
% Since also $\vec{x}^{(1)},\vec{x}^{(2)}\in{}\hat{X}_0$, we now have
% $\vec{x}^{(1)}=\vec{x}^{(2)}$.  
It is unclear whether it is possible, for more general $A_r$, to have
more than one fixed point.

Now define $r_0\ge0$ to be the maximum value of $r\in\RR$ such that
$\sum_{r\le{}r_0}A_r\kappa_r\le{}C$.  Suppose that the reservation
parameters $k_1,\dots,k_r$ are allowed to increase.  Further
consideration of the processes~$\pi_\vec{x}$ shows that, in the limit
(formally as these reservation parameters tend to infinity), the fixed
point $\hat{\vec{x}}$ is necessarily unique and is such that
$P_r(\hat{\vec{x}})=1$ for all $r\le{}r_0$, with, in the heavy traffic
case, $0\le{}P_{r_0+1}(\hat{\vec{x}})\le1$ and $P_r(\hat{\vec{x}})=0$
for all $r\ge{}r_0+2$.  Since the stationary distributions
associated with our sequence of networks converge to that concentrated
on the unique fixed point~$\hat{\vec{x}}$, it follows that the
reservation strategy does indeed approximate, and in the limit
achieve, the complete prioritisation of call types discussed in
Section~\ref{sec:simple-model}.  As mentioned there, and as easily
verified from the above analysis, quite small values of the
reservation parameters $k_1,\dots,k_r$ are sufficient to achieve a
very good approximation to this prioritisation.

Even in the present single-resource case it is possible to achieve
nonuniqueness of the fixed points of the limit
process~$\vec{x}(\cdot)$ by the use of more general, and sufficiently
perverse, controls, in particular with the use of acceptance sets of
the form $\fc{\AA}_r=\{\fc{m}\colon{}A_r\le\fc{m}\le{}k_r+A_r\}$ for
some $k_r\ge0$ (see Bean \textit{et al.}, 1997).  Thus we may
construct networks which
have several (very different) regimes which are quasi-stationary in
the sense discussed above.

\subsection{Multi-resource networks: the uncontrolled case}
\label{sec:multi-reso-netw}

We now consider multi-resource networks, and again study the behaviour
of the fluid limit process~$\vec{x}(\cdot)$ associated with the Kelly
limiting scheme.  Here in general a rich variety of behaviour is
possible.  However, in the case of the uncontrolled networks of
Section~\ref{sec:rev}, in which calls of all types are accepted
subject only to the availability of sufficient capacity, the process
$\vec{x}(\cdot)$ is rather well-behaved.  Note that here, in terms of
the available free capacity, the acceptance sets are given by, for
each $r\in\RR$,
\begin{equation}
  \label{eq:35}
  \fc{\AA}_r = \{\fc{\vec{m}}\colon\ \fc{m}_j\ge A_{jr} \text{ for all $j$}\}.
\end{equation}

Recall also that $X$ is as given by \eqref{eq:30}.  Define the
(real-valued) concave function $f$ on $X$ by
\begin{equation}
  \label{eq:36}
  f(\vec{x}) = \sum_{r\in\RR}(x_r\log\nu_r - x_r\log\mu_rx_r + x_r)
\end{equation}
and let $\hat{\vec{x}}$ be the value of $\vec{x}$ which maximises
$f(\vec{x})$ in $X$.  Kelly (1986) shows that, as
$N\rightarrow\infty$, the stationary distribution of the
process~$\vec{x}^N(\cdot)$ converges to that concentrated on the
single point~$\hat{\vec{x}}$.  (Indeed this is the basis of his
original derivation of the the limiting acceptance probabilities
considered in Section~\ref{sec:appr-large-netw}.)

Assume for the moment the unique existence of the functions $P_r$ on
$X$ introduced above.  Then, for the fluid limit
process~$\vec{x}(\cdot)$, it follows from \eqref{eq:31} and
\eqref{eq:36} that $df(\vec{x}(t))/dt=g(\vec{x}(t))$ where the
function $g$ on $X$ is given by
\begin{align*}
   g(\vec{x})  
  & = \sum_{r\in\RR}\frac{\partial f(\vec{x})}{\partial x_r}%
  \Bigl(\nu_rP_r(\vec{x}) - \mu_r x_r\Bigr)\\
  & = \sum_{r\in\RR}\Bigl(\log \nu_r - \log\mu_r x_r\Bigr)
  \Bigl(\nu_rP_r(\vec{x}) - \mu_r x_r\Bigr).
\end{align*}
Analogously to the preceding section, for each $\vec{x}\in{}X$, the
limiting acceptance probabilities $P_r(\vec{x})$ are given by
consideration of the stationary distribution of a ``free capacity''
Markov process whose transition rates depend on $\vec{x}$.  Some
simple analysis of the equilibrium equations which define this
stationary distribution (see Zachary, 2000) now shows that
$g(\vec{x})\ge0$ for all $\vec{x}\in{}X$ with equality if and only if
$\vec{x}=\hat{\vec{x}}$.

Thus the dynamics of the limit process~$\vec{x}(\cdot)$ are such that,
away from the point~$\hat{\vec{x}}$, the function $f(\vec{x}(\cdot))$
is always strictly increasing.  It thus acts as a Lyapunov function,
ensuring that all trajectories of the process~$\vec{x}(\cdot)$
converge to the single fixed point~$\hat{\vec{x}}$.  Indeed a rigorous
application of the fluid limit theory of Hunt and Kurtz (1994) (again
see Zachary, 2000, for details) shows this result continues to hold
even if the functions~$P_r$ on $X$ are not uniquely defined (whether
this can ever happen in the case of uncontrolled networks remains an
open problem).  The result therefore establishes an important
stability property of uncontrolled networks, and guarantees that the
stationary distribution describes the typical behaviour of the
network.

\subsection{Multi-resource networks: the general case}
\label{sec:multi-reso-netw-1}

For general multi-resource networks, the fluid limit
process~$\vec{x}(\cdot)$ associated with the Kelly limiting scheme may
fail to be unique, and may in particular exhibit multiple fixed
points.  We describe in some detail an elementary example, which is a
simplification of one due to Hunt (1995b).  Suppose that $R=3$, $J=2$,
and that the matrix~$A$ is given by
\begin{displaymath}
  A =
  \begin{pmatrix}
    1 & 0 & 1\\
    0 & 1 & 1
  \end{pmatrix}
  .
\end{displaymath}
Thus in particular calls of types $1$ and $2$ each require capacity
from a single resource, while calls of type~$3$ require capacity from
both resources in the network.  Suppose further that the (free
capacity) acceptance
sets are given by, for some $k_1,k_2\ge1$,
\begin{displaymath}
  \fc{\AA}_1 = \{\fc{\vec{m}}\colon 1\le \fc{m}_1\le k_1\},
  \quad
  \fc{\AA}_2 = \{\fc{\vec{m}}\colon 1\le \fc{m}_2\le k_2\},
  \quad
  \fc{\AA}_3 = \{\fc{\vec{m}}\colon \fc{m}_1\ge1,\,\fc{m}_2\ge 1\}.
\end{displaymath}
(As Hunt remarks, this is not entirely unrealistic: in more complex
networks, operating under some form of alternative routing, certain
resources may have calls of certain types routed over them precisely
when the network is in general very busy.)  Finally suppose that
$\mu_r=1$ for all $r$ and that the vectors~$\vec{\nu}$ and $\vec{C}$
defined in Section~\ref{sec:fluid-limits-large} (each to be scaled by
$N$ for the $N$th member of the sequence of networks) are given by
$\vec{\nu}=(\nu_1,\nu_2,\nu_3)$ and $\vec{C}=(C,C)$.

The process~$\vec{x}(\cdot)$ takes values in the space
$X=\{\vec{x}\in\Rp^3\colon x_1+x_3\le C,~x_2+x_3\le C\}$.  Its
dynamics may be determined through the fluid limit theory outlined
above.  For
$\vec{x}\in{}X_0:=\{\vec{x}\in{}X\colon{}x_1+x_3<C,~x_2+x_3<C\}$
(corresponding to limit points of the dynamics well away from the
capacity constraints) the limiting acceptance probabilities are
well-defined and given by
\begin{equation}
  \label{eq:38}
  P_1(\vec{x})= P_2(\vec{x})=0,\qquad P_3(\vec{x})=1.
\end{equation}
For $\vec{x}\in{}X_1:=\{\vec{x}\in{}X\colon{}x_1+x_3=C,~x_2+x_3<C\}$
and for
$\vec{x}\in{}X_2:=\{\vec{x}\in{}X\colon{}x_1+x_3<C,~x_2+x_3=C\}$
(corresponding in both cases to limit points of the dynamics such that
only one capacity constraint is relevant) the limiting acceptance
probabilities are again well-defined and given by consideration of a
Markov process on $\Zp$ as in the single resource case considered in
Section~\ref{sec:single-reso-netw}.  (For $\vec{x}\in{}X_1$, for
example, it follows from the definition of $\fc{\AA}_2$ that the
transition rates of this Markov process are as if $\nu_2=0$.)  For
$\vec{x}\in{}X_{12}:=\{\vec{x}\in{}X\colon{}x_1+x_3=C,~x_2+x_3=C\}$ it
is necessary to consider also a ``free capacity'' Markov process on
$\Zp^2$.

In the case $\nu_3\le{}C$, these Markov processes all fail to possess
stationary distributions and the limiting acceptance probabilities are
given by \eqref{eq:38} for all $\vec{x}\in{}X$.  Thus the limit
process~$\vec{x}(\cdot)$ is as if $\nu_1=\nu_2=0$ and all trajectories
of this process are deterministic functions of their initial values and
tend to the single fixed point $\hat{\vec{x}}=(0,0,\nu_3)$.

The case $\nu_3>C$ is more interesting.  Here it is readily verified
that the limit process~$\vec{x}(\cdot)$ possesses no fixed points in
$X_0$.  Within $X_1$ consideration of the stationary distribution of
the Markov process defined in Section~\ref{sec:single-reso-netw} shows
that there is a single fixed point $\vec{x}^{(1)}=(a_1,0,C-a_1)$ for
some $a_1$ which is independent of $\nu_2$.  Similarly within $X_2$
there is a single fixed point $\vec{x}^{(2)}=(0,a_2,C-a_2)$ for some
$a_2$ which is independent of $\nu_1$.  However, within $X_{12}$ the
dynamics of the limit process~$\vec{x}(\cdot)$ are not deterministic.
It is further not difficult to show that all trajectories of
$\vec{x}(\cdot)$ which avoid the set~$X_{12}$ tend deterministically
to one of the two fixed points~$\vec{x}^{(1)}$, $\vec{x}^{(2)}$ above
(depending on whether the set $X_1$ or the set $X_2$ is hit first).
Those trajectories of $\vec{x}(\cdot)$ which do hit $X_{12}$ may, in
an appropriate probabilistic sense, tend to either $\vec{x}^{(1)}$ or
$\vec{x}^{(2)}$.

The interpretation of the above behaviour is the following.  Suppose
that $N$ is large and that, for example, resource~$1$ fills to
capacity first.  Then this resource remains full and blocks sufficient
of the type~$3$ calls to ensure that resource~$2$ remains only
partially utilised, with no calls of type~$2$ ever being accepted.
This corresponds to a ``quasi-stationary'' state whose limit, as
$N\to\infty$, is concentrated on the fixed point~$\vec{x}^{(1)}$.
Alternatively, if resource~$2$ fills to capacity first, the network
settles, for an extended period of time, to a quasi-stationary state
whose limit is concentrated on the fixed point~$\vec{x}^{(2)}$.
While, for finite $N$, transitions between these two quasi-stationary
states will eventually occur, the time taken to do so can be shown to
increase exponentially in $N$.

The behaviour in the above example is typical of that which may occur
in more general networks---in particular those using alternative
routing strategies---which are poorly controlled.  Fluid limits may be
used to study behaviour in networks with high capacities and
correspondingly high arrival rates, and to choose values of, for
example, reservation parameters so as to ensure that the network does
not spend extended periods of time in states in which it is operating
inefficiently.  A realistic example here is the fully-connected
network with alternative routing considered in
Section~\ref{sec:mult-reso-models}.

As noted above, fluid limits may also be used to study
equilibrium behaviour, especially in the case where all trajectories
of the limit process~$\vec{x}(\cdot)$ tend to a unique fixed
point~$\hat{\vec{x}}$.  In particular we may show that, for the Kelly
limiting regime considered here, the limiting stationary distribution
of the free capacity processes~$\vec{m}^N(\cdot)$ in general only has
a product form in the case of uncontrolled networks.  This
product-form assumption is the basis of the commonly used
approximations considered in Section~\ref{sec:mult-reso-models}.  Its
justification owes more to the results for the diverse routing limit
also considered there and in
Section~\ref{sec:diverse-rout-limit}.

\subsection{The diverse routing limit}
\label{sec:diverse-rout-limit}

In this section we consider the fluid limit obtained under the diverse 
routing regime discussed in the Introduction.  Although a high degree 
of symmetry is required in order to obtain formal limits, the results 
obtained lend support to the commonly made assumptions of independence 
of resource blocking which are used, for example, in the construction 
of the approximations discussed in Section 3.2.

As outlined earlier, the diverse routing regime holds when the numbers 
of resources and possible ``routes'' in the network increase, while
the total capacity
and arrival rate at each resource remains constant.  
%This limit is 
%equivalent to the thermodynamic limit described in, for instance,
%Graham and Meleard (1995).  
For this limit to 
exist we require a high degree of symmetry in the network. There are two 
canonical examples (with variants) that have been extensively studied. 
We describe both here using the terminology of communications networks.

The first is the so-called \emph{star network} (see, for instance,
Whitt, 1985,
Ziedins and Kelly, 1989, Hunt, 1995a). 
Here there are $K$ links, each with capacity $C$. 
The scale parameter of the regime is then taken to be $K$. 
Assume that calls of any {\em size} $r \geq 1$ require unit capacity 
at each of $r$ resources and have holding times with unit mean.
Then in a symmetric network there are $\binom{K}{r}$
possible choices of the set of links for such a call.  
Let the arrival rate for each such choice be 
$\nu_r^{K} = \nu_r/\binom{K-1}{r-1}$,
so that the total arrival rate at each resource for calls of size $r$ is 
exactly $\lambda_r$.  For example, we may assume that the $K$ links 
are distributed around a central hub, through which all communications 
must pass.  Many variants of this model are 
possible---multiple call sizes can coexist in the network, as can multiple 
capacities, provided only that the proportion of links with any given 
capacity remains constant as $K$ increases.  The network is assumed to
have fixed
routing and the only permissible controls are those on admission.

Let $\vec{x}^{K} (t) = (x_j^{K} (t),\,j\in\JJ)$ where $x_j^{K}(t)$ is
the proportion
of links in which $j$ units of capacity are in use at time~$t$. For the 
network without admission controls, Whitt (1985) showed that, given the 
initial point $\vec{x}(0)$, the process~$\vec{x}^{K}(\cdot)$ converges 
weakly to a deterministic limit process $\vec{x}(\cdot)$, which satisfies 
a set of first-order differential equations with a unique fixed point 
$\hat{\vec{x}}$, such that $\vec{x}(t) \rightarrow \hat{\vec{x}}$ as 
$t \rightarrow \infty$ for all initial $\vec{x}(0)$.  The limit 
$\hat{\vec{x}}$ coincides exactly with that given by the Erlang fixed point 
approximation.  Recall that the latter is obtained from the assumption that 
the stationary free capacity distributions on the various links of the
network are independent of each other. For the case where all calls are of 
size two, Hunt (1995a) obtained a functional central limit theorem for the 
process~$\vec{x}^{K}(\cdot)$, with the limit an Ornstein-Uhlenbeck diffusion 
process (as previously conjectured by Whitt), which was then extended
to more general sizes and initial conditions by Graham and Meleard (1995).  
In the case of networks with 
admission controls very little has been proved.  MacPhee and Ziedins (1996) 
studied such networks and gave a weak convergence result for the process 
$\vec{x}^{K}(\cdot)$.  However, there remain many open questions about the 
behaviour of this process.

The second canonical example of the diverse routing regime is that of the 
fully connected network with alternative routing (Hunt and Laws, 1993).  
Here both admission and routing controls are possible.  The network has 
$N$ nodes; between each pair of these there is a link with capacity $C$,
so that the total number of links is $K = \binom{N}{2}$.  Here again $K$ 
is the scale parameter. Calls arrive at each link at rate $\nu$; each call 
has a unit capacity requirement and holding time of mean 1.  There are 
three possible actions on the arrival of a call: (i) accept the call
at that link, (ii) select a
pair of links that form an alternative route between that pair of nodes 
and route the call along this, or (iii) reject the call. Hunt and Laws showed
that an asymptotically optimal policy, in the sense of minimising the average 
number of lost calls in equilibrium, is to route a call directly if possible 
and otherwise to route it via an alternative route, provided that the
remaining
free capacity on each link of the alternative route is at least some
reservation
parameter $k$, where the optimal choice of $k$ is determined by the parameters 
$K$ and $\nu$. The optimal choice of alternative route is given by choosing 
that which is least busy, i.e.\ which maximises the minimum of the free 
capacities on the two links.  The analysis of Hunt and Laws largely dispenses 
with the graph structure inherent in the choice of alternative routes, an 
assumption justified by analogy with earlier results of Crametz and
Hunt (1991) 
in relation to the simpler model without reservation. 

As in the example of the star network, of interest here is the process
$\vec{x}^{K} (\cdot)$, defined as earlier.  Hunt and Laws showed weak 
convergence of this process to a deterministic limit process. 
They showed that this limit process satisfies differential equations 
which yield the constraints for a linear programming problem, the 
solution to which gives an upper bound on the acceptance probabilities. 
(These constraints correspond to the detailed balance equations that
in equilibrium govern the changes in occupancy of a single link.)
They further showed that their policy achieves this upper bound.

\section{Further developments and open questions}
\label{sec:recent-devel-vari}

Our discussion has of necessity omitted many topics of interest, 
some of which we mention briefly here, as well as discussing some 
remaining open questions.

%Other approaches to the control of loss networks exist in the
%literature.  
%Kelly (1988) has investigated the use of implied
%costs to determine optimal controls, usually assuming that the
%factorisation~\eqref{eq:6} holds.  
%There is an extensive literature using Markov
%decision processes to obtain results on optimal controls 
%(see e.g.\ Key, 1990,
%\textsc{add Feinberg and Reiman, 1994 + some others?}).
%Both of these methods produce
%controls that, while tailored to known arrival rates, may not
%be robust to changes in these rates.  Additionally, the
%Markov decision approach is typically computationally very
%expensive, and produces state-dependent policies that are
%usually complicated to implement.

One such topic is the application of large deviations 
techniques to loss networks in order to estimate, for example,
blocking probabilities in cases where it is important to keep these
very small.  For an excellent introduction 
to this see, for instance, Shwartz and Weiss (1994); later 
papers include those by Simonian \textit{et al.} (1997) and by Graham
and O'Connell (2000). 

In some models of communications networks,
particularly those whose graph structure is tree-like, the network
topology may be
such as to lend itself to more accurate calculations of
acceptance probabilities, involving recursions that do not make 
the link independence assumption~\eqref{eq:45} that is such an
essential feature of
the approximations presented above (see Zachary and Ziedins, 1999).

Extensions of loss network models include
recent work by Antunes \textit{et al.} (2005)
which studies a variant of the model where customers may 
obtain service sequentially at a number of resources, each of which 
is a loss system.  The aim here is to model a 
cellular wireless system where a call in progress may move from base 
station to base station.  Several authors have also considered
explicitly systems with time-varying arrival rates and/or retries 
(see, for example, Jennings and Massey, 1997,
and Abdalla and Boucherie, 2002).

A large number of interesting and important open problems remain.  The
approach to most of these seems to lie in a better understanding of
network dynamics.  There has been no systematic investigation of how
to achieve asymptotically optimal control in a general network (for
example in the sense of Section~\ref{sec:mult-reso-models}), using
controls which are simple, decentralised, and robust with respect to
variations in network parameters, although, for communications
networks, there is a belief that this will usually combine some form
of alternative routing with the use of reservation parameters to
guarantee stability.

A further major problem is that of the identification of instability,
where the state of a network may remain over extended periods of
time in each of a number of ``quasi-equilibrium'' distributions, some
of which may correspond to highly inefficient performance.
Instability is further closely linked to problems of phase transition
in the probabilistic models of statistical physics, and to the study
of how phenomena such as congestion propagate through a network.  At
present results only exist for some very regular network topologies
(see, for instance, Ramanan \emph{et al.}, 2002 and Luen \emph{et al.}, 2006).

Questions related to those above concern the identification of fluid
limits, and in particular the problem of the uniqueness of their
trajectories given initial conditions.  It is notable that the
uniqueness question has not yet been resolved even in the case of a
general \emph{uncontrolled} loss network, although it is known that
here all trajectories do tend to the same fixed point, thus
guaranteeing network stability.  Further, while fixed points of fluid
limits identify quasi-equilibrium states of a network, detailed
behaviour within such states, and the estimation of the time taken to
pass between them, requires a more delicate analysis based on the
study of diffusion limits.  Here relatively little work has been done
(see Fricker \emph{et al.}, 2003).
% and Kelly and Williams, 2004).

Finally we mention that loss networks may be seen as a subclass of a
more general class of stochastic models, with state space $\Zp^R$ for
some $R$ and fairly regular transition rates between neighbouring
states.  Notably their analysis has much in common with that of
processor-sharing networks, in which calls again have a simultaneous
resource requirement.  A unified treatment is still awaited.

\end{document}